\documentclass[draft,12pt]{article}
\usepackage{enumerate}
\usepackage{amsthm,amsmath,amsfonts}

\def\ZZ{\mathbb{Z}}
\def\RR{\mathbb{R}}
\def\NN{\mathbb{N}}
\def\vect#1{\mathbf{#1}}
\def\tr{\top }

\theoremstyle{plain}
\newtheorem{proposition}{Proposition}
\newtheorem{theorem}[proposition]{Theorem}
\newtheorem{lemma}[proposition]{Lemma}
\newtheorem{corollary}[proposition]{Corollary}

\theoremstyle{definition}
\newtheorem*{definition}{Definition}

\theoremstyle{definition}

\newtheorem{remark}{Remark}

\begin{document}

\title{On Fuglede's conjecture and the existence of universal spectra}
\author{B\'alint Farkas, M\'at\'e Matolcsi, P\'eter M\'ora}

\maketitle
\begin{abstract}
Recent methods developed by Tao \cite{tao}, Kolountzakis and Matolcsi \cite{nspec} have led to
counterexamples to Fugelde's Spectral Set Conjecture in both directions. Namely, in $\RR^5$ Tao
produced a spectral set which is not a tile, while Kolountzakis and Matolcsi showed an example
of a non-spectral tile. In search of lower dimensional non-spectral tiles we were led to
investigate the Universal Spectrum Conjecture (USC) of Lagarias and Wang \cite{lagwang}. In
particular, we prove here that the USC and the ``tile $\rightarrow$ spectral'' direction of
Fuglede's conjecture are equivalent in any dimensions. Also, we show by an example that the
sufficient condition of Lagarias and Szab\'o \cite{lagszab} for the existence of universal
spectra is not necessary. This fact causes considerable difficulties in producing lower
dimensional examples of tiles which have no spectra. We overcome these difficulties by invoking
some ideas of R\'ev\'esz and Farkas \cite{revfark}, and obtain non-spectral tiles in $\RR^3$.

Fuglede's conjecture and the Universal Spectrum Conjecture remains open in 1 and 2 dimensions.
The 1 dimensional case is closely related to a number theoretical conjecture on tilings by
Coven and Meyerowitz \cite{covmey}.
\end{abstract}



\section{Introduction}

We briefly summarize the appearing concepts and some related known results. Let $G$ denote some
finite Abelian group or $\ZZ^d$ or $\RR^d$. (The notions below can be extended to the more
general setting of locally compact Abelian groups, but we will restrict our attention to the
mentioned cases.) We will always consider the standard Haar measures on these groups. The
cyclic group of $n$ elements will be denoted by $\ZZ_n$. The dimension of a finite Abelian
group $G$ is understood to be the smallest $d$ such that $G$ is a factor group of $\ZZ^d$. The
Fourier zero-set $Z_T$ of a set $T\subset G$ is defined as the zero-set of the Fourier
transform of the indicator function of $T$, i.e., $Z_T:= \{v\in \Hat G: \ \Hat \chi_T (v)=0\}$;
its complement is denoted by $Z_T^c:=\Hat G\setminus Z_T$.

\begin{definition}A set $\Lambda\subseteq\Hat G$ is called a {\it spectrum} of a bounded open set $T\subseteq G$ if the restricted characters $\{ \lambda (z)|_{z\in T} \}_{\lambda\in \Lambda}$
form a complete orthogonal system in $L_2(T)$. A bounded open set  $T\subseteq G$ is called
{\it spectral} if it possesses a spectrum.
\end{definition}
It is immediate from the definition that for a subset $T$ of a finite group $G$, $\Lambda
\subset \Hat G$ is a spectrum of $T$ if and only if $\Lambda -\Lambda \subset Z_T\cup \{0\}$.

The other class of sets under study is that of {\it translational tiles}:

\begin{definition}A bounded open subset $T$ of $G$ is said to be a \emph{tile}
if the whole group $G$ can be covered by translated disjoint copies of $T$ up to a set of zero
measure. That is, there exists a set $T'\subseteq G$, called a \emph{tiling complement of $T$}
such that $T'+T=\cup_{t'\in T'} (t'+T) =G$ up to a set of zero measure, and the union is
assumed to be disjoint, i.e.,  $T-T\cap T'-T'=\{0\}$.
\end{definition}

In a finite group $G$ there are two well-known (easy) necessary and sufficient conditions for
$T'$ to be a tiling complement of $T$: one is that $|T|\cdot|T'|=|G|$ and $T-T\cap
T'-T'=\{0\}$, while the other is that $|T|\cdot|T'|=|G|$ and $Z_T\cup Z_{T'}=\Hat G\setminus
\{0\}$.

Connecting these two notions Fuglede's spectral set conjecture \cite{fug} asserts that a
bounded open subset of $\RR^d$ is a tile if and only if it is spectral.

Fuglede \cite{fug} proved the special case of the conjecture, when the spectrum or the tiling
complement is assumed to be a lattice. Later on, for many years several positive results seemed
to indicate the validity of the conjecture (see e.g. \cite{convexcurv,convexplane,spectralsym,
konyagin,laba1, laba2, lagwang,unispectra}. Recently, however, Tao \cite{tao} disproved the
``spectral $\rightarrow$ tile'' direction in $\RR^5$ and higher dimensions. Kolountzakis and
Matolcsi \cite{hadamard} reduced this dimension to 3, and also constructed a counterexample to
the ``tile $\rightarrow$ spectral'' direction in $\RR^5$, see \cite{nspec}. R\'ev\'esz and
Farkas refined the arguments of \cite{nspec} to produce a non-spectral tile in $\RR^4$. One of
the aims of this paper is to remove the existing discrepancy between the dimensions by
presenting an example of a tile in $\RR^3$ which does not have a spectrum. In doing so, we will
also investigate the relation of Fuglede's conjecture to the Universal Spectrum Conjecture of
Lagarias and Wang \cite{lagwang} and a certain sufficient condition of Lagarias and Szab\'o
\cite{lagszab} on the existence of universal spectra. The arguments are based on combinatorial
and  Fourier analytic conditions of tiling and spectrality but, unfortunately, some numerical
calculations cannot be avoided. The methods developed in this paper may be useful later in 1 or
2 dimensional considerations. In fact, it is known \cite{konyagin}
 that in $\ZZ$ the ``tile $\rightarrow$ spectral'' direction of Fuglede's conjecture would follow form a number theoretical conjecture on tilings by Coven and Meyerowitz \cite{covmey}. Therefore, any example of a non-spectral tile in $\ZZ$ would immediately disprove the Coven--Meyerowitz conjecture, too. However, despite some numerical experiments in 1-dimensional cyclic groups, we have not been able to find such examples yet.

First, we review and complement some relevant results from the literature. The general
approach, developed first by Tao \cite{tao}, is to find a counterexample in a finite Abelian
group and then transfer the example to the lattice $\ZZ^d$ and, finally, to $\RR^d$.  We recall
the following results:

\begin{theorem}[\cite{nspec}]\label{thm:rd}
Suppose $B \subseteq \ZZ^d$ is a finite set and $Q = (0,1)^d$ is the unit cube. Then $B$ is a
spectral set in $\ZZ^d$ if and only if $B+Q$ is a spectral set in $\RR^d$.
\end{theorem}

It is also clear that if $A$ tiles $\ZZ^d$ then $B+Q$ tiles $\RR^d$.
 This means that any counterexample in $\ZZ^d$ to the ``tile $\rightarrow$ spectral'' direction of Fuglede's conjecture will automatically lead to a counterexample in the corresponding Euclidean space $\RR^d$.

\begin{theorem}[\cite{nspec}]\label{thm:zd}
Let $\vect n = (n_1,\ldots,n_d) \in \NN^d$ and consider a set $A \subseteq G =
\ZZ_{n_1}\times\cdots\times\ZZ_{n_d}$. For the set
\begin{equation}\label{ft-grid} T = T(\vect n, k)
= \{0,n_1,2n_1,\ldots,(k-1)n_1\}\times\cdots\times \{0,n_d,2n_d,\ldots,(k-1)n_d\}
\end{equation}
define $B(k) = A + T$. Then, for large enough values of $k$, the set $B(k)\subset\ZZ^d$ is
spectral in $\ZZ^d$ if and only if $A$ is spectral in $G$.
\end{theorem}

It is also clear that if $A$ tiles $G$ then $B(k)$ tiles $\ZZ^d$ for all values of $k$. This
means that any counterexample in a finite group to the ``tile $\rightarrow$ spectral''
direction of Fuglede's conjecture will automatically lead to a counterexample in the
corresponding lattice $\ZZ^d$. Therefore, our task is to construct a non-spectral tile in some
3-dimensional finite group.

\begin{remark} The corresponding results also hold in this generality in the other direction of Fuglede's conjecture. Namely, any spectral set in $\ZZ^d$ which is not a tile leads (by addition of the unit cube) to a spectral set in $\RR^d$ which is not a tile. Also, any spectral set in a finite group $G =
\ZZ_{n_1}\times\cdots\times\ZZ_{n_d}$ which is not a tile leads (by constructing $B(k)$ as
above) to a spectral set in $\ZZ^d$ which is not a tile. The proof of the first statement is
fairly straightforward, while the second is essentially contained in \cite{tao} and also in
\cite{matolcsi}, Proposition 2.1 and 2.5. We will not need these statements here, but it is
instructive to see that in both directions of Fuglede's conjecture any finite group
counterexample can be transferred automatically to the Euclidean setting, where the conjecture
was originally formulated.
\end{remark}

\section{The equivalence of the Universal Spectrum Conjecture and Fuglede's conjecture}

In order to construct a non-spectral tile in a 3-dimensional finite group we will need the
notion of universal spectrum (see \cite{lagwang}). This notion was originally defined in the
Euclidean setting \cite{lagwang}, but all known sufficient conditions for the existence of
universal spectra go back to considerations in finite groups \cite{lagwang, lagszab}.
Therefore, in this paper, for the sake of simplicity we will remain in the finite setting,
which will suffice for our purposes.

\begin{definition} A subset $S\subset\Hat G$ is a {\it universal spectrum} of a set $T\subset G=\ZZ_{n_1}\times\cdots\times\ZZ_{n_d}$ if $S$ is a spectrum of all tiling complements $T'$ of $T$ in $G$.
\end{definition}

The Universal Spectrum Conjecture (USC) of Lagarias and Wang stated that in any finite Abelian
group $G$ all tiles $T$ posses a universal spectrum. This conjecture was first disproved in
\cite{nspec} in $\ZZ_6^5$, and later in \cite{revfark} in $\ZZ_6^4$, as an essential step in
producing non-spectral tiles in $\RR^5$ and $\RR^4$, respectively. We now prove that the USC
and the ``tile $\rightarrow$ spectral'' direction of Fuglede's conjecture are equivalent in the
sense that the failure of one in any dimension will automatically result in the failure of the
other in the same dimension (i.e., it is not by chance that the examples of \cite{nspec} and
\cite{revfark} led to counterexamples to Fuglede's conjecture).

\begin{proposition}\label{equi}
For any dimension $d$, the Universal Spectrum Conjecture is valid for all $d$-dimensional
finite groups if and only if all tiles are spectral sets in all $d$-dimensional finite groups.
\end{proposition}
\begin{proof}
One direction of this statement is trivial. Namely, if $T$ is a non-spectral tile in a group
$G$ then any tiling complement $T'$ does not possess a universal spectrum in $\Hat G$.

Conversely, assume that we find a $d$-dimensional group
$G=\ZZ_{n_1}\times\cdots\times\ZZ_{n_d}$ and a tile $T\subset G$ which does not have a
universal spectrum. Let $k:=|T|$ and $n:=|G|$. We will exhibit a non-spectral tile $R$ in a
larger (but still $d$-dimensional) group $G_1:=\ZZ_{n_1}\times\cdots\times\ZZ_{n_d}\times
\ZZ_p$, where $p$ is a large integer, relatively prime to $n_1, \dots, n_d$.

Note that $S\subset \Hat G$ is a universal spectrum of $T$ if and only if $|S|=n/k$ and
$S-S\subset \cap_{j=1}^m Z_{T_j'}\cup\{0\}$, where $T_j'$ run through all possible tiling
complements of $T$. By assumption $T$ does not have a universal spectrum, which implies that
for any set $S\subset\Hat G$, $|S|=n/k$ we have a ``witness'' $\vect v_S \in S-S$ such that
$\vect v_S\notin \cap_j Z_{T_j'}\cup\{0\}$. Let $\vect v_1, \vect v_2, \dots, \vect v_r$ denote
the finite set of all such witnesses. Consider now the matrix
$$
A=
\begin{pmatrix}
\Hat\chi_{T_1'}(\vect v_1)&\Hat\chi_{T_2'}(\vect v_1)&\cdots&\Hat\chi_{T_m'}(\vect v_1)\\
\Hat\chi_{T_1'}(\vect v_2)&\Hat\chi_{T_2'}(\vect v_2)&\cdots&\Hat\chi_{T_m'}(\vect v_2)\\
\vdots&&\ddots&\vdots\\
 \Hat\chi_{T_1'}(\vect v_r)&\Hat\chi_{T_2'}(\vect v_r)&\cdots&\Hat\chi_{T_m'}(\vect v_r)
\end{pmatrix}.
$$
We know that each row contains a non-zero entry. We now choose an integer vector $\vect k:=
(k_1, k_2, \dots, k_m)^\tr$ such that $A\vect k\ne 0$ and $k_1+k_2+\cdots+ k_m=p$ is relatively
prime to $n_1, n_2, \dots, n_d$. (It is easy to see that such choice is possible, as the $A{\bf
k}\ne 0$ condition means only an exclusion of $r$ hyperplanes, and the relative prime condition
means only an exclusion of a set of density strictly less than 1.)

We will now glue together the desired non-spectral tile $R\subset G_1$ from several copies of
the sets $T_1', \dots, T_m'$. The idea is that we can consider $G_1$ as $p$ ``layers'' of $G$
and we will copy the sets $T_j'$ on different layers.

 We can regard the elements of $G_1$ as column vectors of length $d+1$. (Note, however, that the
  dimension of $G_1$ is still $d$ as $p$ was chosen relatively prime to $n_1, \dots, n_d$; in fact
  it would suffice that $p$ is relatively prime to one them.) Also, the elements of $\Hat G$ can be
  regarded as row vectors, the action of a character $\gamma\in \hat G$ on an element $\vect x\in G$ being
  defined as $\gamma (\vect x):=
e^{\sum_{j=1}^{d+1} \gamma_j x_j/n_j}$ (where $n_{d+1}:=p$). Let $\vect z_j=(0,0,\dots
,j)^\tr$. For any set $A\subset G$ the notation $\tilde{A}$ will stand for the set $A$ extended
by zero in the last coordinate. Let also $1=\sigma_1\leq\sigma_2\leq\cdots\leq \sigma_p= m$ be
a sequence of integers, the number $i$ occurring exactly $k_i$ times among $\sigma_j$ (recall
that $k_1+k_2+\cdots+k_m=p$). Consider the set
$$
R=\bigcup_{j=1}^p\left(\vect z_j+\tilde{T_{\sigma_j}'}\right)
$$
We claim that $R$ is a tile in $G_1$ and it is not spectral. It is clear that $R$ tiles $G_1$
because a tiling complement can be given as $\tilde{T}$.

Consider any set $L\subset\Hat G_1$, $|L|=|R|$ as a candidate for being a spectrum of $R$. By
the pigeonhole principle there exist an $L_1\subset L$, $|L_1|=n/k$ such that the last
coordinates of the elements of $L_1$ are equal. Consider the set $\tilde{S_1}$ whose elements have the
same coordinates as those of $L_1$ except for the last coordinate which is set to 0 in $\tilde{S_1}$.
Then $\tilde{S_1}-\tilde{S_1}=L_1-L_1\subset L-L$. Consider now the witness $\vect v_{S_1}$ corresponding to
$S_1$, and the extended vector $\tilde{\vect v_{S_1}}$. We have
$$
\Hat\chi_{R}(\tilde{\vect v_{S_1}})=k_1\cdot \Hat\chi_{T_1'}(\vect v_{S_1})+k_2\cdot
\Hat\chi_{T_2'}(\vect v_{S_1})+\cdots+k_m\cdot \Hat\chi_{T_m'}(\vect v_{S_1})\ne 0
$$
by construction. This shows that $R$ is not spectral in $G_1$ and the proof is complete.
\end{proof}
\begin{remark}
One can also introduce the notion of {\it universal tiling complement}. A set $U\subset \Hat G$
is a universal tiling complement of $T\subset G$ if $U$ is a tiling complement in $\Hat G$ of
all spectra of $T$.

Then one can prove the ``dual'' of the statement above, i.e., that all spectral sets are tiles
in all $d$-dimensional finite groups if and only if all spectral sets possess universal tiling
complements in all $d$-dimensional finite groups. In fact, one can use an analogous
construction as above, building up layer by layer a spectral set which is not a tile in a
larger group $G_1=G\times \ZZ_p$ (to see that the constructed set does not tile $G_1$ one needs
to recall the Fourier condition $Z_T\cup Z_{T'}=G\setminus \{0\}$ of tiling pairs). We do not
give a detailed proof here as we will not need this result. However, it is very well possible
that this statement can be useful in producing 1 or 2 dimensional examples in the future.
\end{remark}

\section{A 3-dimensional tile without spectrum}

The results of the previous section show that our task is reduced to finding a tile $T$ of a
3-dimensional finite group $G$ which does not have a universal spectrum. We will exhibit such a
set in $\ZZ_{24}^3$. Unfortunately, it is not at all straightforward to check whether a set $T$
possesses universal spectra or not. There is an elegant sufficient condition by Lagarias and
Szab\'o \cite{lagszab}:

\begin{proposition}\label{sufficient}
For a given set $T$ in a finite group $G$, if a set $S\subset\Hat G$ satisfies the conditions
$|S|=|G|/|T|$ and $S-S\subset Z_T^c$ then $S$ is a universal spectrum of $T$, and also $S$ is a
universal tiling complement of $T$.
\end{proposition}
\begin{proof}
We need to prove the second part of the proposition, which is not contained in \cite{lagszab}.
Assume $L$ is any spectrum of $T$. Then $|L|\cdot|S|=\Hat G$ and $L-L\cap S-S=\{0\}$, because
$L-L\subset Z_T\cup\{0\}$ and $S-S\subset Z_T^c$. It follows that $L+S=\Hat G$.
\end{proof}
In fact, in \cite{lagszab} it is tentatively conjectured that the existence of such set $S$ is
also a necessary condition for the existence of universal spectrum. If it were so, we could
simply use the duality argument of \cite{nspec} to produce a set without universal spectrum in
$G=\ZZ_{24}^3$. The idea is to use the {\it mod 8} log-Hadamard matrix $K$ given in
\cite{hadamard} (i.e., the matrix with entries $e^{2\pi i K_{j,k}}$ is a complex Hadamard
matrix containing 8$^{\text{th}}$ roots of unity):
$$
K:=\frac{1}{8}\left (
\begin{array}{cccccc}
0&0&0&0&0&0\\
0&4&2&6&6&2\\
0&2&4&1&5&6\\
0&6&3&4&2&7\\
0&6&7&2&4&3\\
0&2&6&5&1&4
\end{array}
\right ).
$$
Then, following a decomposition in \cite{hadamard}, one can define a spectral set $T_1$ in
$\ZZ_{24}^3$ with spectrum $L$ as (again the columns are the elements of $G$, while the rows
correspond to elements of the dual group $\Hat G$)
$$
T_1:=\left (
\begin{array}{cccccc}
0&2&4&1&5&6\\
0&6&3&4&2&7\\
0&6&7&2&4&3
\end{array}
\right ) \ \ \ \ \mathrm{and}  \ \ \ L:=3 \left (
\begin{array}{cccccc}
0&0&0\\
0&1&1\\
1&0&0\\
0&1&0\\
0&0&1\\
7&1&1
\end{array}
\right ).
$$
Note that $24K=LT_1$ {\it mod 24}, therefore $L$ is indeed a spectrum of $T_1$ in $G$. Note also,
that $L$ is contained in the subgroup of elements whose coordinates are all divisible by 3.
This subgroup has $8^3$ elements, hence $L$ cannot tile this subgroup due to obvious
divisibility reasons. It is also well-known (and easy) that this implies that $L$ cannot tile
$\Hat G$ either.

It is not hard to see that $T_1$ tiles $G$ (this can be seen e.g. via the homomorphism $\varphi : G\to \ZZ_{24}$ induced by the row vector $(2,9,3)$, but we will need to
modify $T_1$ later anyway), but the existence of a set $S\subset \Hat G$, $|S|=24^3/6$ and
$S-S\subset Z_{T_1}^c$ is impossible due to the following reason: such an $S$ would be a tiling
complement of $L$ by Proposition \ref{sufficient}, which is impossible as $L$ does not tile
$\Hat G$.

If the sufficient condition of Proposition \ref{sufficient} were also necessary then we could
conclude that $T_1$ does not have a universal spectrum in $G$.
 We will show in the Appendix, however, by means of a particular example, that the condition of Proposition
  \ref{sufficient} is not necessary. Of course it still might happen that the set $T_1$ above does not
  have a universal spectrum but, in any case, we are unable to check it at the time of writing.
  (In general, even the elegant sufficient condition of Proposition \ref{sufficient} seems to be
   hard to check algorithmically in large groups, let alone finding all tiling complements of a given set.)
   The failure of the necessity of the Lagarias--Szab\'o condition poses some difficulty in checking whether a
    set possesses universal spectra, and therefore presents an obstacle to finding a 3-dimensional counterexample
    to Fuglede's conjecture. We will use ideas of Farkas and R\'ev\'esz to overcome
this difficulty. The observation is that we are free to add +8 or +16 to the entries of $T_1$
without ruining the decomposition $24K=LT_1$ {\it mod 24}. We must find an alteration $T$ of
$T_1$ such that the existence of a universal spectrum of $T$ can be excluded.

\begin{proposition}\label{usc}
The set
$$
T:=\left (
\begin{array}{cccccc}
0&10&20&1&21&14\\
0&22&3&20&2&7\\
0&22&23&18&4&11
\end{array}
\right )$$ is a tile in $G=\ZZ_{24}^3$ which does not have a universal spectrum.
\end{proposition}
\begin{proof}
As observed before, the decomposition $24K=LT$ still holds, therefore $L$ is a spectrum of $T$.

Consider all the {\it mod 24} vectors $\vect v_{ij}:=\vect l_i-\vect l_j\in \Hat G$ where
$\vect l_i, \vect l_j$ are arbitrary rows of the matrix $L$. For each such vector $\vect
v_{ij}$ we will exhibit a tiling complement $T_{ij}'$ of $T$ in $G$ in such a way that
$v_{ij}\notin Z_{T_{ij}'}$. Accepting the existence of such $T_{ij}'$ for the moment, we can
easily show that $T$ does not have a universal spectrum. Indeed, if $S$ were a universal
spectrum, then $|S|=|G|/|T|$ and $S-S\subset \cap_{ij} Z_{T_{ij}'}\cup\{0\}$ would hold, and
therefore $S-S\cap L-L=\{0\}$ would follow. That is, $S+L$ would be a tiling of $\Hat G$, which
is a contradiction because $L$ is not a tile, as observed in the paragraph after the definition
of $L$. It remains to show the existence of $T_{ij}'$.

Consider all possible {\it mod 8} differences $\vect k_i-\vect k_j$ of the rows of the integer
matrix $8K$. Let $K-K$ denote the matrix containing these differences as row vectors. Now,
regard $K-K$ as a {\it mod 24} matrix and modify the entries by +8 or +16 in such a way that
each row becomes a tile in $\ZZ_{24}$, and also the {\it mod 3} rank of the resulting matrix
$P$ is 3. It will soon be apparent why these modifications are helpful in finding the sets
$T_{ij}'$. We remark that the existence of such modifications appears to be pure luck. We give
a possible example below:

$$K-K=\begin{pmatrix}
    0&0&2&4&4&6 \\
    0&0&4&2&6&4 \\
    0&0&4&6&2&4 \\
    0&0&6&4&4&2 \\
    0&2&1&6&4&5 \\
    0&2&4&1&5&6 \\
    0&2&5&4&6&1 \\
    0&2&6&5&1&4 \\
    0&4&1&5&3&7 \\
    0&4&2&6&6&2 \\
    0&4&3&1&7&5 \\
    0&4&5&7&1&3 \\
    0&4&6&2&2&6 \\
    0&4&7&3&5&1 \\
    0&6&2&3&7&4 \\
    0&6&3&4&2&7 \\
    0&6&4&7&3&2 \\
    0&6&7&2&4&3
\end{pmatrix}\longrightarrow
\begin{pmatrix}
    0&16&2&4&12&14 \\
    0&16&12&2&14&4 \\
    0&16&12&14&2&4 \\
    0&16&14&12&4&2 \\
    0&2&1&14&12&13 \\
    0&2&12&1&13&14 \\
    0&2&13&12&14&1 \\
    0&2&14&13&1&12 \\
    0&12&1&13&11&23 \\
    0&12&2&22&14&10 \\
    0&12&11&1&23&13 \\
    0&12&13&23&1&11 \\
    0&12&22&2&10&14 \\
    0&12&23&11&13&1 \\
    0&22&10&11&23&12 \\
    0&22&11&12&10&23 \\
    0&22&12&23&11&10 \\
    0&22&23&10&12&11
\end{pmatrix}=P$$

It is easy to check that all required properties are fulfilled. In fact, each row of the
modified matrix $P$ has tiling complement $C_1=\{0,3,6,9\}$ or $C_2=\{0,1,6,7\}$ in $\ZZ_{24}$,
and regarding $P$ {\it mod 3} an easy Gaussian elimination shows that the 1$^{\text{st}}$,
2$^{\text{nd}}$ and 4$^{\text{th}}$ rows $\vect p_1, \vect p_2, \vect p_4$ generate the others.

Observe that the set $T$ above is defined in such a way that the rows coincide {\it mod 3} with
$\vect p_1, \vect p_2, \vect p_4$ (and, of course, the entries of $T$ coincide {\it mod 8} with
those of $T_1$).

Consider now an arbitrary row vector $\vect v_{ij}=\vect l_i-\vect l_j$. We will exhibit the
existence of the required tiling complement $T_{ij}'$. For the sake of clarity we follow the
proof through a particular example: let $\vect v_{31}=\vect l_3-\vect
l_1=(3,0,0)-(0,0,0)=(3,0,0)$. Take the corresponding row $\vect k_i-\vect k_j$ of $K-K$, i.e.,
$\vect k_3-\vect k_1=(0,2,4,1,5,6)$ in our particular case. Consider the corresponding row
$\vect p_{ij}$ of the matrix $P$, i.e., $(0,2,12,1,13,14)$ in our case. We claim that there
exists a {\it mod 24} row vector $\vect y_{ij}$ which is a solution of the equation $\vect
y_{ij}T=\vect p_{ij}$ {\it mod 24}. Clearly, a solution of the same equation {\it mod 3}
exists, as $\vect p_{ij}$ is in the linear span of the rows of $T$ {\it mod 3} (recall that $T$
was chosen in such a way that its rows generate {\it every} vector $\vect p_{ij}$ {\it mod 3}).
In our case the {\it mod 3} solution is seen to be $(0,2,0)$. A solution of the same equation
{\it mod 8} is simply obtained by dividing each entry of $\vect v_{ij}$ by 3, i.e., in our case
a {\it mod 8} solution is $(3,0,0)/3=(1,0,0)$. (This is because $\frac{1}{3}LT=8K$ {\it mod
8}.) Then, a solution $\vect y_{ij}$ {\it mod 24} can easily be obtained from the {\it mod 3}
and {\it mod 8} solutions; in our example it is $\vect y_{31}=(9,8,0)$.

Given such $\vect y_{ij}$ we can define a homomorphism $\varphi_{ij}:\ZZ_{24}^3 \to \ZZ_{24}$
by the formula $\varphi_{ij}(\vect x):=\langle \vect y_{ij}, \vect x\rangle$. This homomorphism
takes the set $T$ to the elements of the row $\vect p_{ij}$ by construction, and this resulting
set tiles $\ZZ_{24}$ with complement $C_{ij}:=C_1$ or $C_{ij}:=C_2$ also by construction. In
our example, $\varphi_{31}(T)=(0,2,12,1,13,14)$, which tiles $\ZZ_{24}$ with complement
$C_{31}:=C_1=\{0,3,6,9\}$. Finally, the desired tiling complement $T_{ij}'$ is defined as the
pre-image of $C_{ij}$ under $\varphi_{ij}$. Here we need to invoke an elementary result of
Szegedy \cite{szegedy}.
\begin{lemma}\label{lem:pb}Let $G$ be a finite Abelian group, $T\subseteq G$ and suppose that
there exists a homomorphism $\varphi:G\to H$ such that $\varphi$ is injective on $T$ and
$\varphi(T)$ is a tile in $H$. Then $T$ tiles also $G$, and a tiling complement is given by
$\varphi^{-1}(T')$ where $T'$ is a complement of $\varphi(T)$.
\end{lemma}

Thus, we define $T_{ij}':=\varphi^{-1}(C_{ij})$. It remains to check that $\vect v_{ij}\notin
Z_{T_{ij}'}$. The point of the whole construction above is that we can now evaluate $\Hat
\chi_{T_{ij}'}(\vect v_{ij})$. Note that each homomorphism $\varphi_{ij}$ is easily seen to be
surjective (indeed, a homomorphism $\varphi(x,y,z):=\langle (a,b,c), (x,y,z)\rangle$ is
\emph{not} surjective if and only if $a,b,c$ are all even or all are divisible by $3$; whereas
our vectors are not of this type). Therefore every element in $\ZZ_{24}$ has $24^2$ pre-images
in $\ZZ_{24}^3$. Observe that $3\vect y_{ij}=\vect v_{ij}$ {\it mod 24}, hence for any $\vect
x\in T_{ij}'$ we have $\langle \vect v_{ij},\vect x \rangle \in 3C_{ij}$. Let $\rho
=(1+i)/\sqrt{2}$ denote the first 8$^{\text{th}}$ root of unity. Then
$$\Hat \chi_{T_{ij}'}(\vect v_{ij})=\sum_{\vect x\in T_{ij}'}e^{2\pi i/24 \langle \vect v_{ij}, \vect x \rangle }=
 \sum_{\vect x\in T_{ij}'}e^{2\pi i/24 \langle 3\vect y_{ij}, \vect x \rangle }= \sum_{\vect x\in T_{ij}'}e^{2\pi i/8 \langle \vect y_{ij}, \vect x \rangle }= 24^2\sum_{k\in C_{ij}}\rho^k \ne 0.$$
The last sum is non-zero as $\rho^0+\rho^3+\rho^6+\rho^9 \ne 0$ and
$\rho^0+\rho^1+\rho^6+\rho^7 \ne 0$.
\end{proof}
\noindent Putting together Propositions \ref{usc} and \ref{equi}, and Theorems \ref{thm:zd} and
\ref{thm:rd} we obtain a 3-dimensional counterexample to Fuglede's ``tile $\rightarrow$
spectral'' conjecture:

\begin{corollary}
There exists an appropriate finite union of unit cubes in $\RR^3$ which tiles the space but
which is not spectral.
\end{corollary}

\begin{remark}
At present, all known counterexamples to Fuglede's conjecture (in either direction, and in any
dimensions) have their origins in the existence of complex Hadamard matrices with certain
properties. It is conceivable that a tile having no universal spectrum (or a spectral set
having no universal tiling complement) can be exhibited in a 1 or 2 dimensional finite group
without any reference to Hadamard matrices. By the results of this paper such an example would
immediately lead to a counterexample to (the corresponding direction of) Fuglede's conjecture.
The 1-dimensional case seems particularly interesting, as it is related to the number
theoretical conjecture of Coven and Meyerowitz. In search of a tile without universal spectrum
we have conducted some numerical experiments in several cyclic groups. The main difficulty is
the lack of quick algorithms for deciding whether a set is a tile, and whether it has universal
spectrum. Given the lack of such algorithms we were unable to search large groups exhaustively,
but our ``sporadic'' tests indicate that such examples, if they exist at all, are to be found
in cyclic groups of fairly large order.
\end{remark}

\section{Appendix}

Here we exhibit an example which shows that the sufficient condition given in Proposition
\ref{sufficient} is not necessary for the existence of a universal spectrum. This example was
found earlier in \cite{revfark}, but there the authors could not decide whether the set given
below has a universal spectrum.

\begin{proposition}
In the group $G=\ZZ_6^4$ the set given by the columns
$$
T:=\left (
\begin{array}{cccccc}
0&1&0&0&0&2\\
0&0&1&0&0&2\\
0&0&0&1&0&2\\
0&0&0&0&1&2
\end{array}
\right )
$$
does have a universal spectrum $U$, while there exists no $S\subset\Hat G$ such that $|S|=6^3$
and $S-S\subset Z_T^c$.
\end{proposition}

To see that such $S$ cannot exist we use the argument that has appeared several times above
(and, in fact, was the basis of the arguments in \cite{nspec}). Observe that $T$ is spectral in
$G$ with a possible spectrum
$$L:=2 \left (
\begin{array}{cccc}
0&0&0&0\\
0&1&1&2\\
1&0&2&2\\
1&2&0&1\\
2&2&1&0\\
2&1&2&1
\end{array}
\right ).
$$

Note that $L$ is contained in the subgroup of elements with even coordinates. This subgroup has
$3^4$ elements, therefore $L$ does not tile this subgroup, and hence does not tile $\Hat G$
either. Now, the existence of a proposed set $S$ would mean, by Proposition \ref{sufficient},
that $S$ is a universal tiling complement of $T$, which is a contradiction as $L$ does not tile
$\Hat G$.

To see that $T$ does have a universal spectrum $U$ is more tricky. In fact, $U$ can be defined as
$U:=\{ (u_1, u_2, u_3, u_4): u_1+u_2+u_3+u_4=0  \ mod \ 6\}$. Then $U-U=U$ and therefore we
only need to show that $U\subset \cap_j Z_{T_{j}'}\cup\{0\}$ where $T_{j}'$ runs through all
the tiling complements of $T$. Write $U=U_0\cup U_1$, where $U_0=U\cap Z_T$ and $U_1=U\cap
Z_T^c$. It is trivial from the Fourier condition $Z_T\cup Z_{T_j'}=G\setminus \{0\}$ that
$U_1\subset \cap_j Z_{T_{j}'}\cup\{0\}$. It is an easy calculation to show that $U_0$ consists
of all coordinate permutations of the vector $\vect v=(2,2,4,4)$. By symmetry of $T$ it is enough to
show that $\vect v\in  \cap_j Z_{T_{j}'}$. This, however, is non-trivial and, unfortunately, we do
not have a neat ``structural'' proof of this fact. In fact, the first proof that we found was
listing out all tiling complements $T_j'$ by a computer search (which, in itself, is a nearly
impossible task due to the large size of the group $\ZZ_6^4$, and the lack of known quick algorithms.)
Below we present a proof that is easy (but tedious) to check by hand.

Consider any tiling complement $T'$ of $T$. We need to show that $\vect v\in Z_{T'}$. Let $P$
denote the following set:
$$P=\bigl\{ (3,4,4,4)^\tr, (4,3,4,4)^\tr, (4,4,3,4)^\tr, (4,4,4,3)^\tr\bigr\}.$$

\noindent We will show the following:

\vskip1ex\noindent\textbf{Fact 1.} If $\vect t\in T'$ then there exists an $\vect x\in P$ such
that $\vect t+ \vect x\in T'$.

\vskip1ex\noindent\textbf{Fact 2.} If $\vect t\in T'$, $\vect x\in P$ and $\vect t+\vect x\in
T'$, then $\vect t+2\vect x\in T'$.

\vskip1ex\noindent\textbf{Fact 3.} If $\vect t\in T'$, $\vect x, \vect y\in P$ and $\vect
t+\vect x\in T'$ and $\vect t+\vect y\in T'$ then $\vect x=\vect y$.

\vskip2ex Assuming these facts for the moment, we can conclude our argument easily. Indeed, the
above statements imply that $T'$ is a disjoint union of 36 6-cycles of the form $C_j=\{\vect
t_j, \vect t_j+ \vect x_j, \vect t_j+2\vect x_j, \vect t_j+3\vect x_j, \vect t_j+4\vect x_j,
\vect t_j+5\vect x_j\}$ where $\vect x_j\in P$. Therefore,
$$\Hat \chi_{T'} (\vect v)=\sum_{j=1}^{36} (e^{2\pi i/6\langle \vect v,\vect t_j\rangle }\sum_{k=0}^5 e^{2\pi i/6\langle \vect v, k\vect x_j\rangle })=0,$$
because the inner sums are all easily seen to be zero for all $x_j\in P$.

\vskip2ex
To show Fact 1 is easy. We can assume without loss of generality that $\vect t=\vect 0$. Let us now see how the tiling $T+T'=\ZZ_6^4$ covers the point $(4,4,4,4)^\tr$. In order to cover $(4,4,4,4)^\tr$, the set $T'$ must contain one of the following points: \\
$\{ (4,4,4,4)^\tr, (3,4,4,4)^\tr, (4,3,4,4)^\tr, (4,4,3,4)^\tr, (4,4,4,3)^\tr, (2,2,2,2)^\tr \}$. However, \\
$(4,4,4,4)^\tr\in T'$ would mean that $T+T'$ covers $\vect 0$ twice, while $(2,2,2,2)^\tr\in
T'$ would mean that $T+T'$ covers $(2,2,2,2)^\tr$ twice.

\vskip2ex Fact 3 is also straightforward by contradiction. We can assume once again that $\vect
t=\vect 0$ and, without loss of generality, that $\vect x=(3,4,4,4)^\tr$, $\vect
y=(4,3,4,4)^\tr$. Then, $\vect t+\vect x=\vect x\in T'$ and $\vect t+\vect y=\vect y\in T'$
would imply that $T+T'$ covers $(4,4,4,4)^\tr$ twice, a contradiction.

\vskip2ex The proof of Fact 2 is generated by a computer algorithm, but can also be checked by
hand. We argue by contradiction. Assume that $\vect t=\vect 0\in T'$, $\vect x=\vect t+\vect
x=(3,4,4,4)^\tr\in T'$, and $2\vect x=\vect t+2\vect x=(0,2,2,2)^\tr\notin T'$. Then
$(1,2,2,2)^\tr$ must be covered by $T+T'$, and there are three ways to do it:
\begin{enumerate}[ii)]
\item[1.] $(1,1,2,2)^\tr\in T'$
\item[2.] $(1,2,1,2)^\tr\in T'$
\item[3.] $(1,2,2,1)^\tr\in T'$
\end{enumerate}
By symmetry of all appearing sets with respect to the last three coordinates, it is enough to
check one of these cases, say 1; i.e. assume that  $(1,1,2,2)^\tr\in T'$. Then $(0,0,0,5)^\tr$
must be covered by $T+T'$, and there are three ways to do it.
\begin{enumerate}[ii)]
\item[1.a.] $(0,5,0,5)^\tr\in T'$
\item[1.b.] $(0,0,5,5)^\tr\in T'$
\item[1.c.] $(0,0,0,4)^\tr\in T'$
\end{enumerate}

\noindent All of these cases lead to contradiction in the following manner:

\vskip1ex\noindent\textbf{Case 1.a.}
           \\   $(0,0,5,0)^\tr$ can only be covered if $(0,0,4,0)^\tr\in T'$, and then
           \\   $(2,2,1,2)^\tr$ can only be covered if $(2,2,1,1)^\tr\in T'$, and then
           \\   $(4,4,3,4)^\tr$ can only be covered if $(4,3,3,4)^\tr\in T'$, and then
           \\   $(2,1,1,2)^\tr$ can only be covered if $(2,0,1,2)^\tr\in T'$, and then
           \\   $(1,5,5,0)^\tr$ can only be covered if $(1,5,5,0)^\tr\in T'$, and then
           \\   $(0,5,5,1)^\tr$ can only be covered if $(5,5,5,1)^\tr\in T'$, and then
           \\   $(5,0,0,1)^\tr$ can only be covered if $(4,0,0,1)^\tr\in T'$, and then
           \\   $(3,1,0,2)^\tr$ can only be covered if $(3,1,5,2)^\tr\in T'$, and then
           \\   $(2,5,4,0)^\tr$ can only be covered if $(0,3,2,4)^\tr\in T'$, and then
           \\   $(4,1,0,2)^\tr$ cannot be covered without covering some point twice.

\vskip1ex\noindent\textbf{Case 1.b.}       \\   $(0,5,0,0)^\tr$ can only be covered if
$(0,4,0,0)^\tr\in T'$, and then
           \\   $(2,1,1,2)^\tr$ can only be covered if $(2,1,0,2)^\tr\in T'$, and then
           \\   $(2,2,1,2)^\tr$ cannot be covered without covering some point twice.

\vskip1ex\noindent\textbf{Case 1.c.}
           \\   $(2,2,2,1)^\tr$ can only be covered if $(2,2,1,1)^\tr\in T'$, and then
           \\   $(1,2,2,1)^\tr$ can only be covered if $(0,2,2,1)^\tr\in T'$, and then
           \\   $(1,2,1,2)^\tr$ can only be covered if $(1,2,0,2)^\tr\in T'$, and then
           \\   $(4,4,4,3)^\tr$ can only be covered if $(4,3,4,3)^\tr\in T'$, and then
           \\   $(5,0,0,5)^\tr$ can only be covered if $(5,0,5,5)^\tr\in T'$, and then
           \\   $(1,3,1,1)^\tr$ can only be covered if $(1,3,1,0)^\tr\in T'$, and then
           \\   $(5,1,0,5)^\tr$ can only be covered if $(5,1,0,5)^\tr\in T'$, and then
           \\   $(5,0,1,5)^\tr$ can only be covered if $(5,5,1,5)^\tr\in T'$, and then
           \\   $(0,5,1,0)^\tr$ can only be covered if $(4,3,5,4)^\tr\in T'$, and then
           \\   $(0,0,1,5)^\tr$ cannot be covered without covering some point twice.

\end{document}